\documentclass[12pt]{amsart}

\usepackage{amsfonts}
\usepackage{amssymb}
\usepackage[T2A]{fontenc}
\usepackage[cp1251]{inputenc}
\usepackage[russian,english]{babel}

\newtheorem{thm}{Theorem}[section]
\newtheorem{lem}[thm]{Lemma}
\newtheorem{cor}[thm]{Corollary}
\newtheorem{defn}[thm]{Definition}
\newtheorem{prop}[thm]{Proposition}
\theoremstyle{definition}

\def\<{\langle}
\def\>{\rangle}

\def\a{\alpha}

\def\e{\varepsilon}
\def\f{{\varphi}}
\def\F{{\Phi}}
\def\A{{\mathcal A}}
\def\M{{\mathcal M}}
\def\cN{{\mathcal N}}

\def\N{{\mathbb N}}
\def\C{{\mathbb C}}
\def\Hom{\mathop{\rm Hom}\nolimits}
\def\bK{{\mathbf{K}}}
\def\bBK{{\mathbf{BK}}}

\begin{document}

\title[Locally unital $C^*$-algebras and seminorm criteria]{Locally unital $C^*$-algebras and limitations of seminorm criteria for compact operators}

\author{D.V. Fufaev}
\email{denis.fufaev@math.msu.ru}
\email{fufaevdv@rambler.ru}

\address{Moscow Center for Fundamental and Applied Mathematics,
\newline
Dept. of Mech. and Math., Lomonosov Moscow State University}

\begin{abstract}
We study the problem of geometric characterization of Banach-compact and $\mathcal A$-compact operators on Hilbert $C^*$-modules in terms of uniform structures generated by seminorms. 
Previously, such a criterion for $\mathcal A$-compactness of adjointable operators was obtained by Troitsky for countably generated range modules with a suitable system of seminorms.
Using locally unital $C^*$-algebras as modules over themselves, we establish fundamental limitations of such characterizations. 
Namely, we show that for nonunital locally unital $C^*$-algebras not only does Troitsky's criterion fail to extend, but, more generally, no system of seminorms can characterize $\mathcal A$-precompactness by total boundedness, nor Banach-compactness of bounded operators or $\mathcal A$-compactness of adjointable operators by total boundedness of the images of the unit ball.
Moreover, the same statement holds for Hilbert $C^*$-modules $\mathcal N$ whenever the algebra $\mathbf K(\mathcal N)$ of $\mathcal A$-compact operators is locally unital and nonunital. We also study Troitsky-type systems of seminorms and their relation to $\mathcal A$-precompactness. We prove that every $\mathcal A$-precompact subset is totally bounded with respect to the strongest of such structures. Consequently, for Hilbert $C^*$-modules admitting a standard frame, all the corresponding notions of total boundedness coincide.\end{abstract}

\keywords{locally unital $C^*$-algebra, Hilbert $C^*$-module, $\mathcal A$-precompactness, total boundedness, Banach-compact operator, $\mathcal A$-compact operator}

\maketitle

\section{Introduction}

The study of uniform structures (or, more precisely, systems of seminorms) on Hilbert $C^*$-modules aimed at characterizing $\A$-compact operators, that is, operators which can be approximated in norm by finite-rank operators, was initiated in the works of Ke\v{c}ki\'c and Lazovi\'c \cite{Lazovic2018, KeckicLazovic2018,KeckicLazovic2019}.
A substantial advance was made by Troitsky \cite{Troitsky2020JMAA}, who introduced a particular uniform structure and obtained a criterion for $\A$-compactness of adjointable operators in terms of this structure. Namely, if $F:\M\to\cN$ is an adjointable operator and $\cN$ is countably generated, then $F$ is $\A$-compact if and only if the image of the unit ball of $\M$ is totally bounded with respect to this uniform structure.

Later, in \cite{TroitFuf2020}, the forward implication was extended to arbitrary Hilbert \(C^*\)-modules, while the converse was established under the assumption that the range module \(\mathcal N\) has a standard frame. This assumption cannot in general be omitted: as shown in \cite{Fuf2021faa,Fuf2022Path}, there exist Hilbert \(C^*\)-modules whose identity operator is not \(\mathcal A\)-compact, while the unit ball is totally bounded with respect to Troitsky's uniform structure (or Troitsky's topology, as in \cite{KeckicLazovic2026}).

To investigate the limitations of such seminorm characterizations we can turn our attention to locally unital \(C^*\)-algebras introduced and studied in \cite{Fuf2025}, where it was proved that a nonunital locally unital \(C^*\)-algebra, regarded as a Hilbert module over itself, admits no frames. 
In the present paper we prove that, for every such algebra $\mathcal A$, the converse implication in Troitsky's criterion fails even for the strongest uniform structure considered here: the unit ball is $(\mathcal A,\mathcal A)'$-totally bounded, while the identity operator is not $\mathcal A$-compact. Moreover, there is no system of seminorms on $\mathcal A$ which 
characterizes \(\mathcal A\)-precompactness by total boundedness, nor Banach-compactness or \(\mathcal A\)-compactness of operators by total boundedness of the images of their unit balls.
The same conclusions hold for Hilbert $C^*$-modules $\cN$ over an arbitrary $C^*$-algebra $\A$ whenever $\mathbf K(\mathcal N)$ is locally unital and nonunital.
Thus, a geometric characterization of $\mathcal A$-compact and Banach-compact operators cannot, in general, be obtained in terms of total boundedness with respect to a system of seminorms.

A closely related direction concerns different variants of these systems of seminorms. In \cite{TroitFufAFA}, the original construction based on admissible systems of module elements was compared with stronger constructions defined by adjointable and by arbitrary bounded \(\mathcal A\)-linear functionals. In particular, total boundedness for the original structure was shown to be equivalent to total boundedness for the structure defined by adjointable functionals.

We strengthen this comparison by proving that every \(\mathcal A\)-precompact subset of an arbitrary Hilbert \(C^*\)-module is \((\mathcal N,\mathcal N)'\)-totally bounded, i.e. totally bounded with respect to the strongest uniform structure considered here. So, for modules admitting a standard frame, all the corresponding notions of total boundedness coincide. 
As a consequence, we obtain a strengthening of the corresponding compactness result: the image of the unit ball under an \(\mathcal A\)-compact operator, and even under a Banach-compact operator, is totally bounded with respect to the strengthened uniform structure.

The paper is organized as follows. In Section~2 we recall the basic necessary facts about Hilbert $C^*$-modules, locally unital $C^*$-algebras, and the uniform structures under consideration.
In Section~3 we prove that every $\mathcal A$-precompact subset of a Hilbert $C^*$-module is $(\mathcal N,\mathcal N)'$-totally bounded and derive the corresponding consequences for modules with standard frames and for Banach-compact operators. 
In Section~4 we study modules for which $\mathbf K(\mathcal N)$ is locally unital and prove the nonexistence results for seminorm characterizations of $\mathcal A$-precompactness, Banach-compactness, and $\mathcal A$-compactness.

\section{Preliminaries}

The basic theory on Hilbert $C^*$-modules can be found in \cite{Lance},\cite{MTBook},\cite{ManuilovTroit2000JMS}. Let us recall some important notions.

\begin{defn}
\rm
Let $\A$ be a $C^*$-algebra.
A pre-Hilbert $C^*$-module over $\A$ is a complex vector space $\M$ which is a (right) $\A$-module (with compatible scalar multiplication: $\lambda(xa)=(\lambda x)a=x(\lambda a)$ for $x\in\M$, $a\in\A$, $\lambda\in\mathbb C$) equipped with an $\A$-valued inner product $\<.,.\>:\M\times\M\to \A$ which is a sesquilinear form such that, for any $x,y\in\M$, $a\in\A$:
\begin{enumerate}
\item $\<x,x\> \ge 0$ and $\<x,x\> = 0$ if and only if $x=0$;
\item $\<y,x\>=\<x,y\>^*$;
\item $\<x,y\cdot a\>=\<x,y\>a$.
\end{enumerate}

This inner product induces the norm $\|x\|=\|\<x,x\>\|^{1/2}$, and a module is called a \emph{Hilbert $C^*$-module} if it is complete with respect to this norm.

Any $C^*$-algebra $\A$ can be considered as a module over itself  with the inner product $\<a,b\>=a^*b$.

\rm
By an \emph{operator} we mean a bounded $\A$-linear map.
An operator is called adjointable if it admits an adjoint (see \cite[Section 2.1]{MTBook}).
We also define the \emph{dual module} $\M'=\Hom_\A(\M,\A)$ as the set of all bounded $\A$-linear functionals (see also \cite[\S2.5]{MTBook}). 
\end{defn}

\begin{defn}\label{dfn:Acompact}
An \emph{elementary} \emph{Banach-compact} operator
$\theta_{x,f}:\M\to\cN$, where $x\in\cN$ 
and $f\in\M'$, is defined as $\theta_{x,f}(z):=x\cdot f(z)$.
Then the Banach space $\bBK(\M,\cN)$
of \emph{Banach-compact} operators
is the norm closure of the subspace generated by all 
elementary Banach-compact operators in the space
$\Hom_\A(\M,\cN)$ of all bounded $\A$-linear maps.
If we take only functionals of the form $\<y,\cdot\>$ for some $y\in\M$, then we obtain the definitions of an elementary $\A$-compact operator $\theta_{x,y}$ and the corresponding Banach space $\bK(\M,\cN)$ 
of $\A$-compact operators.
\end{defn}

Among the equivalent definitions of $\A$-precompactness established in \cite{TroitFuf26JMAA}, we shall use the following one. A bounded subset $Z\subset\cN$ is $\A$-precompact if, for every $\e>0$, there exist $z_1,\dots,z_n,y_1,\dots,y_n\in\cN$ such that
$$
\sup_{x\in Z}
||x-\sum_{j=1}^n z_j\langle y_j,x\rangle||<\varepsilon.
$$

This property yields the following Grothendieck-type characterization of the $\A$-compact and Banach-compact operators.

\begin{thm}\label{thm:opercomp} (\cite[Theorem 4.8]{TroitFuf26JMAA})
Let $F: \M\to \cN$ be a bounded (resp., adjointable) $\A$-linear operator and let $B$ be the unit ball of $\M$.
Then the following are equivalent:
\begin{enumerate}
\item[\rm (a)] $F$ is Banach-compact (resp. $\A$-compact).
\item[\rm (b)] $F(B)$ is $\A$-precompact in $\cN$.
\end{enumerate}
\end{thm}

The notion of frames in Hilbert $C^*$-modules was introduced by Frank and Larson in \cite{FrankLarson1999}, \cite{FrankLarson2002}:
\begin{defn}\label{dfn:fr}
Let $\cN$ be a Hilbert $C^*$-module over a $C^*$-algebra $\A$ and $J$ be some index set. A family $\{x_j\}_{j\in J}$ of elements of $\cN$ is said to be a \emph{standard frame} in $\cN$ if there exist positive constants $c_1, c_2$ such that for any $x\in\cN$ the series $\sum\limits_j \<x,x_j\>\<x_j,x\>$
converges in norm in $\A$ (in particular, only countably many terms of the series are nonzero), and the following inequalities hold:
$$c_1\<x,x\>\le \sum\limits_j \<x,x_j\>\<x_j,x\>\le c_2\<x,x\>.$$
\end{defn}

We now recall some facts about local units and locally unital $C^*$-algebras that were introduced in \cite{Fuf2025}, where several examples were also described.
Among them are the commutative algebra $C_0([0,\omega_1))$ of continuous functions vanishing at infinity on the space of all at most countable ordinals, and the algebra $B_\omega(H)$ of bounded operators with separable range on a nonseparable Hilbert space $H$. 

\begin{defn}
A self-adjoint element $h$ in a $C^*$-algebra $\A$ is called a local unit for an element $a\in\A$ if $ah=ha=a$. 
\end{defn}

By \cite[Proposition 2.2]{Fuf2025} we can always assume that for a local unit $h$ the inequalities  $0\le h \le1$ hold.

\begin{defn}
We will say that a $C^*$-algebra $\A$ is locally unital if every $a\in\A$ has a local unit. 
\end{defn}

In \cite[Theorem 3.3]{Fuf2025} a topological characterization of such algebras in the commutative case was obtained. Namely, a $C^*$-algebra $C_0(K)$ is locally unital if and only if the closure of every $\sigma$-compact subset of $K$ is compact (here $K$ is a locally compact Hausdorff space).

From \cite[Proposition 2.4]{Fuf2025} it follows that in such algebras there exists a common local unit for every at most countable set of positive elements. Now we prove that this holds for arbitrary elements.

\begin{lem}
    Let $\A$ be a $C^*$-algebra and $a,h\in\A$. 
    \begin{enumerate}
        \item If $haa^*=aa^*$, then $ha=a$.
        \item If $a^*ah=a^*a$, then $ah=a$.
    \end{enumerate}
    \end{lem}
\begin{proof}
1) Suppose that $haa^*=aa^*$, i.e. $(1-h)aa^*=0$ in the unitalization of $\A$. Then 

$$
||(1-h)a||^2=||(1-h)aa^*(1-h^*)||=0, 
$$
so $(1-h)a=0$, i.e. $ha=a$. The point 2) is similar.
\end{proof}

\begin{prop}\label{countable-common-local-unit}
Let $\A$ be a locally unital $C^*$-algebra. Then for any countable family of its elements $\{a_j\}_{j\in\N}$ there exists their common local unit $h\in\A$.
\end{prop}

\begin{proof}
    An element $h\in\A$ is a local unit for $b\in\A$, $b\ne0$, if and only if it is a local unit for $\frac{b}{||b||}$, so we can assume that $||a_j||\le1$, $j\in\N$.

    Consider $a=\sum\limits_{j=1}^\infty\frac{1}{2^j}(a_j^*a_j+a_ja_j^*)\in\A$ and choose a local unit $h\in\A$ for $a$. By \cite[Proposition 2.4]{Fuf2025} $h$ is a local unit for all $a_ja_j^*$ and $a_j^*a_j$, $j\in\N$; in particular, $a_j^*a_jh=a_j^*a_j$ and $a_ja_j^*=ha_ja_j^*$. By the previous lemma, $ha_j=a_j=a_jh$, i.e. $h$ is a local unit for all $a_j$, $j\in\N$.
\end{proof}

Now introduce the most general uniform structure considered in this paper, defined by bounded \(\mathcal A\)-linear functionals (see \cite{TroitFufAFA} for details).

\begin{defn}\label{dfn:admissyst*}
Let $\cN$ be a Hilbert $C^*$-module over $\A$. A countable system $F=\{f_{i}\}\subset\cN'$ is called $(\cN,\cN^0)'$-admissible for a submodule $\cN^0\subseteq  \cN$ (or $*$-$\cN^0$-admissible, as denoted in \cite{TroitFufAFA}) if
\begin{enumerate}
\item[1)]  for each $x\in\cN^0$ the partial sums of the series $\sum_i (f_i(x))^* f_i(x)$ 
are bounded by $\<x,x\>$;
\item[2)] this series is norm convergent; 
\item[3)]  $\|f_i\|\le 1$ for any $i$. 
\end{enumerate}
\end{defn}

Denote by $\F$ a sequence $\{\f_1,\f_2,\dots\}$ of states on $\A$. For each pair $(F,\F)$ with an $(\cN,\cN^0)'$-admissible $F$, consider the following seminorm

\begin{equation}\label{eq:defnnorm*}
\nu_{F,\F}(x)^2:=\sup_k 
\sum_{i=k}^\infty |\f_k\left(f_i (x)\right)|^2,\quad x\in \cN^0 
\end{equation}
and the corresponding pseudometric $d_{F,\F}(x,y):=\nu_{F,\F}(x-y)$.
In particular, it was proved in \cite{TroitFufAFA} that $\nu_{F,\F}(\cdot)\leq \|\cdot\|$.

\begin{defn}\label{dfn:totbaundset*}
A set $Y\subseteq \cN^0 \subseteq  \cN$ is $(\cN,\cN^0)'$-\emph{totally bounded} if for any $(F,\F)$, where $F \subseteq  \cN'$ is $(\cN,\cN^0)'$-admissible,
it is $\nu_{F,\F}$-totally bounded, i.e. for any $\e>0$ there exists a finite $\e$-net in $Y$ for $\nu_{F,\F}$. The uniform structure generated by the pseudometrics $d_{F,\F}$ will be called the $(\cN,\cN^0)'$-uniform structure.
\end{defn}

Now we describe two important special cases of this construction. The first is obtained by restricting the admissible systems to adjointable functionals, 
i.e., to those \(\mathcal A\)-linear operators \(f:\mathcal N\to\mathcal A\) which admit an adjoint.

Finally, Troitsky's original construction (from \cite{Troitsky2020JMAA}) can be obtained by using the functionals represented by elements of \(\mathcal N\).
Indeed, $\cN$ is included in $\cN'$ by the formula $\widehat{x}(z)=\<x,z\>$ and we have $\<x,z\>=\<z,x\>^*$ and $\f(a^*)=\overline{\f(a)}$ for any $a\in\A$ and any state $\f$ on $\A$. Hence, for $\{x_i\}=X\subset\cN$ we have $|\f(\<x,x_i\>)|=|\f(\<x_i,x\>^*)|=|\f(\<x_i,x\>)|=|\f(\widehat{x}_i(x))|$, where $\widehat{x}_i\in\cN'$. 

The formula (\ref{eq:defnnorm*}) in this case takes the form

\begin{equation}\label{eq:defnnorm}
\nu_{X,\F}(x)^2:=\sup_k 
\sum_{i=k}^\infty |\f_k\left(\<x,x_i\>\right)|^2,\quad x\in \cN^0. 
\end{equation}

If, in the definitions of admissible systems and total boundedness, we restrict ourselves to adjointable \(\mathcal A\)-linear functionals, we will say that the system is \((\mathcal N,\mathcal N^0)'_{ad}\)-admissible and that the set is \((\mathcal N,\mathcal N^0)'_{ad}\)-totally bounded. If we restrict further to functionals represented by elements of \(\mathcal N\) itself, we will say that the system is \((\mathcal N,\mathcal N^0)\)-admissible and that the set is \((\mathcal N,\mathcal N^0)\)-totally bounded.

In \cite[Theorem 2.2]{TroitFufAFA} it was shown that \((\mathcal N,\mathcal N^0)\)-total boundedness and \((\mathcal N,\mathcal N^0)'_{ad}\)-total boundedness coincide.

The following is a special case of \cite[Proposition 1 and Theorem 5]{BanGabPr2016}, but we briefly give a proof in terms of seminorms.

\begin{lem}\label{lem_compseq}
Let $X$ be a vector space, $\nu$ a seminorm on $X$ and $Z\subset X$ be a subset. Then $Z$ is $\nu$-totally bounded if and only if every countable subset of $Z$ is $\nu$-totally bounded.
\end{lem}

\begin{proof}
The forward implication is obvious, let us prove the converse. Suppose $Z$ is not $\nu$-totally bounded, then for some $\e>0$ there exists a countable set $\{z_j\}_{j\in\N}\subset Z$ such that $\nu(z_j-z_i)\ge\e$ for $j\ne i$, in particular, it is not $\nu$-totally bounded. Hence, if every countable subset of $Z$ is $\nu$-totally bounded, then $Z$ is $\nu$-totally bounded too.
\end{proof}

\section{\texorpdfstring{$\A$-precompactness implies $(\cN,\cN)'$-total boundedness}{A-precompactness implies total boundedness}}

One can compare the notions of $(\cN,\cN^0)$-, $(\cN,\cN^0)'_{ad}$- and $(\cN,\cN^0)'$-total boundedness. As mentioned above, by \cite{TroitFufAFA}, $(\cN,\cN^0)$- and $(\cN,\cN^0)'_{ad}$-total boundedness coincide. Also note that if $\cN^0\subset\cN^1\subset\cN^2$, then $(\cN^2,\cN^0)$-total boundedness is stronger than $(\cN^1,\cN^0)$-total boundedness. But in this case 
each $f\in(\cN^2)'$ can be restricted to $\cN^1$,
so $(\cN^1,\cN^0)'$-total boundedness is stronger than $(\cN^2,\cN^0)'$-total boundedness. So, among these conditions $(\cN^0,\cN^0)$-total boundedness is the weakest and $(\cN^0,\cN^0)'$-total boundedness is the strongest. 

It turns out that when $\cN^0$ has a standard frame $(\cN^0,\cN^0)$-total boundedness implies $(\cN^0,\cN^0)'$-total boundedness, that is, they coincide. To simplify the notation we will write $\cN$ in place of $\cN^0$.

\begin{lem}\label{lem_main}
Let $\cN$ be a Hilbert $C^*$-module over a $C^*$-algebra $\A$, $F=\{f_i\}\subset\cN'$ an arbitrary $(\cN,\cN)'$-admissible system and $\F=\{\f_k\}$ an arbitrary sequence of states. Suppose that $Z\subset \cN$ is a bounded set such that for any $\e>0$ there exists $N\in\N$ with

$$
\sup\limits_{x,y\in Z}\sup\limits_{k\in\N}\sum\limits_{i=N+1}^\infty|\f_k(f_i(x-y))|^2\le\e^2,
$$
then $Z$ is $\nu_{F,\F}$-totally bounded.
\end{lem}

\begin{proof}
Assume $Z$ is not $\nu_{F,\F}$-totally bounded; in particular, there exist $\varepsilon>0$ and a sequence $\{b_r\}_{r\in\mathbb N}\subset Z$ 
such that $\nu_{F,\F}(b_r-b_s)\ge\varepsilon$ for $r\ne s$. Define $b_{r,s}=b_r-b_s$.

We can choose $N\in\N$ such that 

$$
\sup\limits_{r,s\in \N}\sup\limits_{k\in\N}\sum\limits_{i=N+1}^\infty|\f_k(f_i(b_r-b_s))|^2\le\frac{\e^2}{4}.
$$

For any $r,s\in\N$, $r\ne s$ we have (noting that
$\sum\limits_{i=k}^{N}(\cdots)=0$ for $k>N$ and that
$\sup\limits_{k\in\N}\sum\limits_{i=k}^N|\f_k(f_i(b_{r,s}))|^2=\max\limits_{k=1,\dots,N}\sum\limits_{i=k}^N|\f_k(f_i(b_{r,s}))|^2$)

$$
\e^2\le\sup\limits_{k\in\N}\sum\limits_{i=k}^\infty|\f_k(f_i(b_{r,s}))|^2\le
\sup\limits_{k\in\N}\left(\sum\limits_{i=k}^N|\f_k(f_i(b_{r,s}))|^2+
\sum\limits_{i=N+1}^\infty|\f_k(f_i(b_{r,s}))|^2\right)\le
$$

$$\le
\sup\limits_{k\in\N}\sum\limits_{i=k}^N|\f_k(f_i(b_{r,s}))|^2+
\sup\limits_{k\in\N}\sum\limits_{i=N+1}^\infty|\f_k(f_i(b_{r,s}))|^2\le
\sup\limits_{k\in\N}\sum\limits_{i=k}^N|\f_k(f_i(b_{r,s}))|^2+\frac{\e^2}{4},
$$
so
$$
\sup\limits_{k\in\N}\sum\limits_{i=k}^N|\f_k(f_i(b_{r,s}))|^2\ge\frac{3}{4}\e^2.
$$

Consider  the finite-di\-men\-si\-on\-al complex space $\C^d$
of dimension $d=N+(N-1)+...+1$ with the norm
\begin{equation}\label{eq:norm_comple}
\|(u_{1,1},\dots,u_{1,N}, u_{2,2},\dots, u_{2,N},\dots, u_{N,N})\|^2
=\max_{1\le k\le N}\sum_{i=k}^N|u_{k,i}|^2
\end{equation}

and the $\C$-linear map 
$$
T: \cN\to \C^d,\qquad y \mapsto \left(\f_k\left(f_i(y) \right)\right),\qquad
k=1,\dots,N,\quad i=k,\dots,N.
$$ 
It is a bounded map because $|\f_k\left(f_i(y) \right)|\le||\f_k||\cdot||f_i||\cdot||y||$, i.e. each coordinate is a bounded map. Thus, $\{T(b_{r})\}_{r\in\N}\subset \C^d$ is a bounded set. On the other hand, for $r\ne s$

$$
||T(b_r)-T(b_s)||^2=||T(b_{r,s})||^2=
\sup\limits_{k\in\N}\sum\limits_{i=k}^N|\f_k(f_i(b_{r,s}))|^2\ge\frac{3}{4}\e^2,
$$
i.e. elements of $\{T(b_r)\}_{r\in\N}$ are $\frac{\sqrt{3}}{2}\e$-separated, which is impossible for a bounded subset of a finite-dimensional space. A contradiction, so $Z$ is $\nu_{F,\F}$-totally bounded.

\end{proof}

\begin{thm}\label{th_acomp_NNtotb}
If $\cN$ is a Hilbert $C^*$-module over a $C^*$-algebra $\A$ and $Z\subset\cN$ is an $\A$-precompact subset, then $Z$ is $(\cN,\cN)'$-totally bounded.
\end{thm}

\begin{proof}
Since $Z$ is $\A$-precompact it is bounded by some constant $M$.
Suppose that $Z$ is not $(\cN,\cN)'$-totally bounded,
i.e. for some $(\cN,\cN)'$-admissible system $F=\{f_i\}_{i\in\mathbb N}$ and some sequence of states $\F=\{\f_k\}_{k\in\mathbb N}$ the set $Z$ is not $\nu_{F,\F}$-totally bounded. In particular, there exist $\delta>0$ and a sequence $\{b_r\}_{r\in\mathbb N}\subset Z$, 
such that $\nu_{F,\F}(b_r-b_s)\ge2\delta$ for $r\ne s$. Define $b_{r,s}=b_r-b_s$ and note $||b_{r,s}||\le2M$.

By $\A$-precompactness 
there exists 
$\{z_k,y_k\}_{k=1}^n$,
$z_k,y_k\in\cN$,
such that
\[
\sup_{x\in Z}\|x-\sum\limits_{k=1}^nz_k\<y_k,x\>\|<\delta/2.
\]

Denote 
$w_{r}=\sum\limits_{k=1}^nz_k\<y_k,b_{r}\>$,
$w_{r,s}=w_{r}-w_{s}=\sum\limits_{k=1}^nz_k\<y_k,b_{r,s}\>$ and observe 

$$
\nu_{F,\F}(w_r-w_s)\ge 
\nu_{F,\F}(b_r-b_s)-\nu_{F,\F}(b_r-w_r)-\nu_{F,\F}(w_s-b_s)\ge
$$

$$\ge
\nu_{F,\F}(b_r-b_s)-||b_r-w_r||-||w_s-b_s||\ge\delta
$$
because $\nu_{F,\F}(x)\le||x||$ for any $x\in\cN$ (as proved in \cite{TroitFufAFA}). In particular, the set $Z_0=\{w_r\}_{r\in\N}$ is also not $\nu_{F,\F}$-totally bounded. But it is bounded in norm because $||w_r-b_r||<\frac{\delta}{2}$ and $||b_r||\le M$, so $||w_r||\le M+\frac{\delta}{2}$.

Since $F=\{f_i\}$ is admissible, for any $z_k$, $k=1,\dots,n$, the series $\sum\limits_{i=1}^\infty(f_i(z_k))^*f_i(z_k)$ converges in norm and for any $\e>0$ we can find $N\in\N$ such that 
$$||\sum\limits_{i=p}^q(f_i(z_k))^*f_i(z_k)||\le\frac{\e^2}{16n^2(M+1)^2\max\limits_{k=1,\dots,n}\{||y_k||+1\}^2}$$
for any $k=1,\dots,n$ and any $q\ge p\ge N+1$ (we add a constant equal to one to avoid division by zero in the possible case $M=0$ or $y_k=0$ for all $k$).

We need the following inequality which holds for any elements $c_1,\dots,c_n$ in any $C^*$-algebra:

$$
\left(\sum\limits_{k=1}^nc_k\right)^*\left(\sum\limits_{k=1}^nc_k\right)\le
n\sum\limits_{k=1}^nc_k^*c_k.
$$
Indeed,

$$
n\sum\limits_{k=1}^nc_k^*c_k-\left(\sum\limits_{k=1}^nc_k\right)^*\left(\sum\limits_{k=1}^nc_k\right)=
\sum\limits_{1\le k<l\le n}(c_k-c_l)^*(c_k-c_l)\ge0.
$$

For any $i,r,s\in\N$ apply this inequality 
with $c_k=f_i(z_k)\<y_k,b_{r,s}\>$,

$$
(f_i(w_{r,s}))^*f_i(w_{r,s})=
\left(\sum\limits_{k=1}^nf_i(z_k)\<y_k,b_{r,s}\>\right)^*\left(\sum\limits_{k=1}^nf_i(z_k)\<y_k,b_{r,s}\>\right)\le
$$

$$\le
n\sum\limits_{k=1}^n\<b_{r,s},y_k\>(f_i(z_k))^*f_i(z_k)\<y_k,b_{r,s}\>
$$
therefore for any $p,q\in\N$, $q\ge p$,

$$
\sum\limits_{i=p}^q(f_i(w_{r,s}))^*f_i(w_{r,s})\le 
n\sum\limits_{k=1}^n\<b_{r,s},y_k\>\sum\limits_{i=p}^q(f_i(z_k))^*f_i(z_k)\<y_k,b_{r,s}\>
$$
and for $q\ge p\ge N+1$

$$
||\sum\limits_{i=p}^q(f_i(w_{r,s}))^*f_i(w_{r,s})||\le
||n\sum\limits_{k=1}^n\<b_{r,s},y_k\>\sum\limits_{i=p}^q(f_i(z_k))^*f_i(z_k)\<y_k,b_{r,s}\>||
\le
$$

$$\le
4nM^2\max_{k=1,\dots,n}\{||y_k||\}^2\sum\limits_{k=1}^n||\sum\limits_{i=p}^q(f_i(z_k))^*f_i(z_k)||
\le \frac{\e^2}{4}.$$

By the Cauchy--Schwarz inequality for states (see, for example, \cite[2.3.11]{BratRob}), for any $k,i,r,s\in\N$ we have

$$
|\f_k(f_i(w_{r,s}))|^2
\le\f_k\left((f_i(w_{r,s}))^*f_i(w_{r,s})\right).
$$

Therefore, for any $p,q\ge N+1$, $q\ge p$, and any $k,r,s\in\N$ we have

$$
\sum\limits_{i=p}^q|\f_k(f_i(w_{r,s}))|^2\le
\sum\limits_{i=p}^q\f_k\left((f_i(w_{r,s}))^*f_i(w_{r,s})\right)=
\f_k\left(\sum\limits_{i=p}^q(f_i(w_{r,s}))^*f_i(w_{r,s})\right)\le
$$

$$\le
||\sum\limits_{i=p}^q(f_i(w_{r,s}))^*f_i(w_{r,s})||\le\frac{\e^2}{4},
$$
hence

$$
\sum\limits_{i=N+1}^\infty|\f_k(f_i(w_{r,s}))|^2\le\frac{\e^2}{4}.
$$

Since $k$ is arbitrary we have

$$
\sup\limits_{k\in\N}\sum\limits_{i=N+1}^\infty|\f_k(f_i(w_{r,s}))|^2\le\frac{\e^2}{4}
$$
for any $r,s\in\N$.
Therefore, we can apply Lemma \ref{lem_main} to the set $Z_0=\{w_r\}_{r\in\N}$ and conclude that it is $\nu_{F,\F}$-totally bounded, a contradiction.
So $Z$ is $\nu_{F,\F}$-totally bounded for any $(F,\F)$ with $(\cN,\cN)'$-admissible $F$, i.e. it is $(\cN,\cN)'$-totally bounded.

\end{proof}

The following result follows from \cite[Theorem 4.7]{TroitFuf26JMAA}.

\begin{thm}
If a Hilbert $C^*$-module $\cN$ over a $C^*$-algebra $\A$ has a standard frame and $Z\subset\cN$ is an $(\cN,\cN)$-totally bounded subset, then $Z$ is $\A$-precompact.
\end{thm}

\begin{cor}
If a Hilbert $C^*$-module $\cN$ over a $C^*$-algebra $\A$ has a standard frame and $Z\subset\cN$ is an $(\cN,\cN)$-totally bounded subset, then $Z$ is $(\cN,\cN)'$-totally bounded.
\end{cor}

By \cite[Theorem 4.8]{TroitFuf26JMAA} the image of the unit ball under a Banach-compact operator is $\A$-precompact.
Hence we obtain the following strengthening of the forward implication in \cite[Theorem 2.13]{Troitsky2020JMAA}.

\begin{cor}
  Let $F: \M\to \cN$ be a bounded $\A$-linear operator. If $F$ is Banach-compact, then $F(B)$ is $(\cN,\cN)'$-totally bounded, where $B$ is the unit ball of $\M$.
\end{cor}

\section{\texorpdfstring{Nonexistence of seminorm criteria for $\A$-precompactness}{Nonexistence of seminorm criteria for A-precompactness}}

In this section we show that, in general, neither \(\mathcal A\)-precompactness of subsets nor Banach-compactness nor \(\mathcal A\)-compactness of operators can be characterized in terms of total boundedness with respect to a system of seminorms. More precisely, we show that there exist Hilbert \(C^*\)-modules for which no system of seminorms can provide such a characterization. In fact, any nonunital locally unital \(C^*\)-algebra \(\mathcal A\), regarded as a Hilbert module over itself, provides such an example: for every system of seminorms on \(\mathcal A\), the corresponding equivalence between \(\mathcal A\)-precompactness and total boundedness fails; similarly, Banach-compactness of bounded operators and \(\mathcal A\)-compactness of adjointable operators cannot be characterized by total boundedness of the images of their unit balls.

Moreover, this is also true for an arbitrary Hilbert $C^*$-module $\cN$ over an arbitrary $C^*$-algebra $\A$ such that the algebra of $\A$-compact operators $\bK(\cN):=\bK(\cN,\cN)$ on $\cN$ is locally unital, but not unital. The case of an algebra considered above can be obtained by taking $\cN=\A$ since $\bK(\A)\cong\A$.

The next two results require only local unitality of \(\mathbf K(\mathcal N)\) (i.e. it may be unital); nonunitality will be needed only afterwards.

\begin{lem}\label{lem:acomp_common_unit}
Suppose that \(\mathbf K(\cN)\) is locally unital. Then for every at most countable subset \(Z=\{\xi_j\}_{j\in\N}\subset\cN\) there exists an operator $h\in\mathbf K(\cN)$ such that
$h\xi=\xi$ for every $\xi\in Z$.
\end{lem}

\begin{proof} By Proposition \ref{countable-common-local-unit} the countable set
$\{\theta_{\xi_j,\xi_j}\}_{j\in\N}
\subset\mathbf K(\cN)$
has a common self-adjoint local unit \(h\in\mathbf K(\cN)\). Hence, for any $j\in\N$,
$$
h\theta_{\xi_j,\xi_j}=\theta_{\xi_j,\xi_j}=\theta_{\xi_j,\xi_j}h.
$$

Denote $\operatorname{Id}-h=g$. Then

$$
g\xi_j\<g\xi_j,\cdot\>=
\theta_{g\xi_j,g\xi_j}=
(\operatorname{Id}-h)\theta_{\xi_j,\xi_j}(\operatorname{Id}-h)=0.
$$
So,

$$
0=\<g\xi_j, \theta_{g\xi_j,g\xi_j} g\xi_j\>=
\<g\xi_j,g\xi_j\<g\xi_j,g\xi_j\>\>=
\<g\xi_j,g\xi_j\>^2=0,
$$
hence, since $\<g\xi_j,g\xi_j\>$ is positive, it is equal to zero, therefore $g\xi_j=(\operatorname{Id}-h)\xi_j=0$, that is, $h\xi_j=\xi_j$.

\end{proof}

\begin{thm}\label{th_balltotb}
Let $\cN$ be a Hilbert $C^*$-module over a $C^*$-algebra $\A$.
Suppose that $\bK(\cN)$ is a locally unital $C^*$-algebra. 

\begin{enumerate}
\item Any countable bounded subset $Z=\{\xi_j\}_{j\in\N}$ of $\cN$ is $\A$-precompact;
\item If $\nu$ is a seminorm on $\cN$ such that every $\A$-precompact subset $Z$ is $\nu$-totally bounded, then 
every bounded subset $Z$ of $\cN$ 
is $\nu$-totally bounded;
\item Every bounded subset of $\cN$ 
is $(\cN,\cN)'$-totally bounded.
\item For any Hilbert $C^*$-module $\M$ over $\A$ and any bounded $\A$-linear operator $F:\M\to\cN$ the image $F(B)$ of the unit ball  $B\subset \M$ is $(\cN,\cN)'$-totally bounded.
\end{enumerate}
\end{thm}

\begin{proof}

(1) By the previous lemma there exists an operator $h\in\bK(\cN)$ such that $h\xi=\xi$ for any $\xi\in Z$.
Since $h$ is $\A$-compact, for any $\e>0$ there exist elements $y_k$, $z_k$, $k=1,\dots,s$, such that 

$$
||h-\sum\limits_{k=1}^s\theta_{z_k,y_k}||<\frac{\e}{\sup\limits_{\xi\in Z}||\xi||+1}.
$$

Using the fact that $h\xi=\xi$ for any $\xi\in Z$ we have

$$
\sup\limits_{\xi\in Z}||\xi-\sum\limits_{k=1}^sz_k\<y_k,\xi\>||=\sup\limits_{\xi\in Z}||(h-\sum\limits_{k=1}^s\theta_{z_k,y_k})\xi||\le||h-\sum\limits_{k=1}^s\theta_{z_k,y_k}||\cdot\sup\limits_{\xi\in Z}||\xi||<\e,
$$
so $Z$ is $\A$-precompact.

(2) By Lemma \ref{lem_compseq} it follows that we need to prove only that any countable subset $Z_0=\{\xi_j\}_{j\in\N}$ of $Z$ is $\nu$-totally bounded.

By item (1) such a countable subset $Z_0$ (which is also bounded since $Z$ is bounded) is $\A$-precompact, so $Z_0$
is $\nu$-totally bounded by the assumption.

(3) follows immediately from (2) and theorem \ref{th_acomp_NNtotb}.

(4) follows from (3) and the fact that for a bounded operator, the image of the unit ball is a bounded subset.
\end{proof}

The following proposition is a corollary of \cite[Theorem 4.8]{TroitFuf26JMAA} and \cite[Theorem 2.4]{Lance}.

\begin{prop}\label{prop:unit-ball-precompact}
For a Hilbert \(\mathcal A\)-module \(\cN\), the following conditions are equivalent:
\begin{enumerate}
\item the unit ball \(B\) is \(\mathcal A\)-precompact;
\item \(\operatorname{Id}\in\mathbf K(\cN)\);
\item the \(C^*\)-algebra \(\mathbf K(\cN)\) is unital.
\end{enumerate}
Consequently, if \(\mathbf K(\cN)\) is nonunital, then \(B\) is not \(\mathcal A\)-precompact.
\end{prop}

The following corollary is a generalization of \cite[Corollary 2.2 and Theorem 2.5]{Fuf2021faa}.

\begin{cor}\label{ne_acomp}
    If $\bK(\cN)$ is nonunital but locally unital, then the identity operator is not $\A$-compact and the unit ball $B$ of $\cN$ is not $\A$-precompact, while it is $(\cN,\cN)'$-totally bounded.
    \end{cor}

Thus, in this situation even the strongest Troitsky-type uniform structure considered here cannot characterize either $\A$-compactness or Banach-compactness by total boundedness of the images of unit balls.

The following theorem shows that the obstruction is not specific to Troitsky's particular choice of seminorms: no alternative system of seminorms can restore such a characterization.

\begin{thm}
Let $\cN$ be a Hilbert $C^*$-module over a $C^*$-algebra $\A$.
Suppose that $\bK(\cN)$ is a locally unital nonunital $C^*$-algebra. 
\begin{enumerate}
\item There doesn't exist a system of seminorms $\{\nu_\a\}_{\a\in\Lambda}$ on $\cN$ such that a bounded subset $Z\subset \cN$ is $\A$-precompact if and only if $Z$ is $\nu_\a$-totally bounded for all $\a\in\Lambda$.
\item There doesn't exist a system of seminorms $\{\nu_\a\}_{\a\in\Lambda}$ on $\cN$ such that a bounded (resp., adjointable) operator $F:\cN\to\cN$ is Banach-compact (resp., $\A$-compact) if and only if $F(B)$ is $\nu_\a$-totally bounded for all $\a\in\Lambda$, where $B$ is the unit ball of $\cN$.
\end{enumerate}
\end{thm}

\begin{proof}
1). Suppose that such a system exists, in particular, if $Z$ is $\A$-precompact, then $Z$ is $\nu_\a$-totally bounded for any $\a\in\Lambda$. Then by theorem \ref{th_balltotb} the unit ball $B$ of $\cN$ is $\nu_\a$-totally bounded for any $\a\in\Lambda$. So, by assumption, $B$ is $\A$-precompact, but it is not true due to Corollary \ref{ne_acomp}. A contradiction, so such a system doesn't exist.

2).
Suppose that such a system exists. Consider for any countable subset $Z_0\subset B$ of the unit ball 
an operator $h\in\bK(\cN)$ such that $h\xi=\xi$ for any $\xi\in Z_0$ (by Lemma \ref{lem:acomp_common_unit}). Then since the image of the unit ball under the $\A$-compact operator $h$ is $\nu_\a$-totally bounded for any $\a\in\Lambda$, we obtain that $Z_0$ is $\nu_\a$-totally bounded for any $\a\in\Lambda$ too. But from lemma \ref{lem_compseq} we obtain that $B=\operatorname{Id}(B)$ is $\nu_\a$-totally bounded for any $\a\in\Lambda$ too, therefore, $\operatorname{Id}:\cN\to\cN$ is at least Banach-compact, but in this case the image of the unit ball $\operatorname{Id}(B)=B$ is $\A$-precompact, which is not true due to corollary \ref{ne_acomp}. 
A contradiction, so such a system doesn't exist.

\end{proof}

We emphasize that the statement above rules out the existence of a system of seminorms yielding a criterion, that is, an equivalence between \(\mathcal A\)-precompactness and total boundedness. Systems of seminorms for which only one of the two implications holds do nevertheless exist. Indeed, the results of the previous section show that the seminorms defining Troitsky-type uniform structures always provide the forward implication: \(\mathcal A\)-precompactness implies \((\mathcal N,\mathcal N)'\)-total boundedness. Conversely, the reverse implication can be obtained, for example, by taking the norm of the Hilbert \(C^*\)-module itself.

If $Z\subset\cN$ is totally bounded with respect to the norm, then for any $\e>0$ there exists a finite set $\{t_k\}_{k=1}^s$ such that for any $x\in Z$ there exists $t_k$ with $||x-t_k||<\e$. By \cite[Proposition 1.3.10]{MTBook} every $t_k$ can be represented as $t_k=z_ka_k$ for some $z_k\in\cN$, $a_k\in\A$. So, $||x-z_ka_k||<\e$ (i.e. all coefficients except the $k$-th are zero) and $||a_k||\le M:=\max\limits_{k=1,\dots,s}||a_k||$, this is the condition a) of \cite[Theorem 4.2]{TroitFuf26JMAA}, hence $Z$ is $\A$-precompact.

\subsection*{Acknowledgements}
The author is grateful to E.V. Troitsky, V.M. Manuilov, A.I. Shtern, K.L. Kozlov, A.Ya. Helemskii and A.I. Korchagin for helpful discussions.

\subsection*{Funding}
This work was supported by the Russian Science Foundation
(Grant No. 24-11-00124).

\end{document}